\documentclass{mscs}
\usepackage{mymacros}

\usepackage{breakcites}

\newtheorem{lemma}{Lemma}[section]
\newtheorem{definition}{Definition}[section]
\newtheorem{theorem}{Theorem}[section]
\newtheorem{proposition}{Proposition}[section]
\newtheorem{remark}{Remark}[section]
\newtheorem{corollary}{Corollary}[section]
\newtheorem{question}{Question}[section]

\title[Synthetic Homology in Homotopy Type Theory]
  {Synthetic Homology in Homotopy Type Theory}
\author[Robert Graham]
  {Robert Graham}
\date{December 2018}
\begin{document}

\global\long\def\help{\operatorname{sgn}}

\global\long\def\colim{\operatorname{colim}}

\global\long\def\merid{\operatorname{merid}}

\global\long\def\smcf{\operatorname{smcf}}

\global\long\def\leftme{\operatorname{left}}

\global\long\def\wglue{\operatorname{wglue}}

\global\long\def\rightme{\operatorname{right}}

\global\long\def\smbase{\operatorname{smbase}}

\global\long\def\smin{\operatorname{smin}}

\global\long\def\smglue{\operatorname{smglue}}

\global\long\def\cfbase{\operatorname{cfbase}}

\global\long\def\cfcod{\operatorname{cfcod}}

\global\long\def\cfglue{\operatorname{cfglue}}

\global\long\def\refl{\operatorname{refl}}

\global\long\def\susp{\operatorname{susp}}

\global\long\def\loopme{\operatorname{loop}}

\global\long\def\merloop{\operatorname{merloop}}

\global\long\def\fib{\operatorname{fib}}

\global\long\def\stab{\operatorname{stab}}

\global\long\def\suspcf{\operatorname{suspcf}}

\global\long\def\suspsm{\operatorname{suspsm}}

\global\long\def\suspsmtwo{\operatorname{suspsm2}}

\global\long\def\susptwo{\operatorname{susp2}}

\global\long\def\smglueleft{\operatorname{smglueleft}}

\global\long\def\smglueright{\operatorname{smglueright}}

\global\long\def\suspthree{\operatorname{susp3}}

\global\long\def\apme{\operatorname{ap}}

\global\long\def\transport{\operatorname{transport}}

\global\long\def\movelr{\operatorname{movelr}}

\global\long\def\pushadj{\operatorname{pushadj}}

\global\long\def\inr{\operatorname{inr}}

\global\long\def\inl{\operatorname{inl}}

\global\long\def\glue{\operatorname{glue}}

\global\long\def\id{\operatorname{id}}

\global\long\def\base{\operatorname{base}}

\global\long\def\cfgluesmbase{\operatorname{cfgluesmbase}}

\global\long\def\cfgluesmin{\operatorname{cfgluesmin}}
\global\long\def\type{\operatorname{type}}
\global\long\def\idtoiso{\operatorname{idtoiso}}
\global\long\def\ua{\operatorname{ua}}

\label{firstpage}
\maketitle
\begin{abstract}
This paper defines homology in homotopy type theory, in the process
stable homotopy groups are also defined. Previous research in synthetic
homotopy theory is relied on, in particular the definition of cohomology.
This work lays the foundation for a computer checked construction
of homology.
\end{abstract}

\tableofcontents
\ifprodtf \newpage \else \vspace*{-1\baselineskip}\fi

%

\section{Introduction}

Homotopy type theory (HoTT) is a new variant of Martin-L$\text{ö}$f's
intuitionistic type theory (ITT) and a potential foundation for all
of mathematics that enables practical computer proof checkers. Key to understanding HoTT is the fact that in ITT,
the identity endows each type with an $\infty$-groupoid structure
\cite{hofmann1998groupoid,lumsdaine2009weak,van2011types}. This
realization prompted the use of homotopy theory techniques to understand
ITT \cite{ voe2010,warren2008homotopy,awodey2009homotopy}. Eventually, this led to the semantics of HoTT \cite{kapulkin2012simplicial}. Formally
HoTT extends ITT with two new concepts: the univalence axiom, and
higher inductive types. The univalence axiom states that any two isomorphic
types are actually equal. Higher inductive types are types in which
the user can specify identities, not just terms. 
 
The univalence axiom and higher inductive types have many different
applications, but of special interest to us is how they relate to
space. Since an $\infty$-groupoid can be viewed as a space, it follows
that types can be viewed as spaces. In this context the univalence
axiom says that any two homotopic spaces are in fact equal, whereas
higher inductive types allow us to define freely generated $\infty$-groupoids.
We can thus do homotopy theory inside of homotopy type theory. This
viewpoint is often called synthetic homotopy theory, as concepts like
paths are primitive instead of being defined in terms of a set of
points (what one could think of as the analytic approach). 

Homotopy type theory enables the ability to use formal proof checkers for theorems that would otherwise be intractable. Initial work in 2013 at the Institute for Advanced Study \cite{hottbook} generated lots of interest  with stunning results like a formalized proof that the fundamental group of the circle is $\mathbb{Z}$ in 100 lines of code. Since then, many other basic homotopy results have been reworked in homotopy type theory, which made homology a natural next step. Specifically, we build off of the following results: the construction of Eilenberg-Maclane spaces \cite{licata2014eilenberg}, showing the torus is the product of two circles \cite{licata2015cubical} (in particular we at times adopt the cubical approach described therein), $\pi_(S^3)\simeq \mathbb{Z}/2\mathbb{Z}$ \cite{brunerie2016homotopy}, and Blakers-Massey theorem \cite{finster2016mechanization}. We especially draw from Evan Cavallo's work on the cohomology case \cite{cavallo2015synthetic} and van Doorn's work on Serre cohomology spectral sequences \cite{van2018formalization}.

We aim to define homology in homotopy type theory, that is,
for each type $X$ we want a group $H_{n}(X)$ that satisfies the
Eilenberg\textendash Steenrod axioms of homology (Defintion \ref{def:homologytheory}).
To this end we would like a definition of $H_{n}$ that does not depend
on the cell structure of a space. Luckily there is such a construction
in the literature, details of which can be found in many introductory texts including Hatcher \cite{hatcher2002algebraic}.
We will broadly follow the approach of Hatcher, however, since our
definitions apply to types, our proofs will often be unique, combining
ideas from type theory with ideas from classical homotopy theory.
As in the case of cohomology, which has already been defined in \cite{cavallo2015synthetic},
we drop the additivity axiom from our definition of homology. However, unlike the case of cohomology, it appears that this should not be necessary. Other researchers have communicated to me progress towards including additivity, but as of now this is not complete to my knowledge.

We have not mechanized these results aside from some basics. Working in HoTT can often involve burdensome computations
involving higher paths. As an example, we originally could not formally prove
certain properties of the smash product that are trivial classically but become very challenging in HoTT (see Lemma \ref{lem:suspsmash}). This issue was ultimately resolved by Floris van Doorn \cite{van2018formalization}. An even greater challenge lies in dealing with colimits, where basic identities of the natural numbers can clog up proofs.

The plan is as follows. To begin with we must first define
stable homotopy groups and prove that they form an homology theory.
To our knowledge stable homotopy theory has not previously been defined
in the literature, but the main tools for working with them have already
been fleshed out (in particular the Freudenthal suspension theorem
(Proposition \ref{prop:(Freudenthal-Suspension-Theorem)}) and the
Blakers-Massey theorem (Theorem \ref{thm:Black Massey})). Our main
result is proving the stable homotopy groups $\pi_{n}^{s}(-)$ form
a homology theory (Corollary \ref{cor:stableishomology}). After which
we extend this to showing $\pi_{n}^{s}(-\land K)$ is a homology theory
for a fixed $K$ (Corollary \ref{cor:stablesm-is-homo}) and finally
to showing $\colim_{i}\pi_{n+i}^{s}(X\land K_{i})$ is a homology
theory for a fixed prespectrum $K_{i}$ (Corollary \ref{cor:colim-is-homo}).
Regular homology $H_{n}(X,G)$ is then defined using the spectrum
of Eilenberg-Maclane spaces $K(G,i)$.

\subsection{Review of Higher Inductive Types}
As mentioned above, higher inductive types (HIT) allow us to define
various spaces in homotopy type theory by specifying the generators
of a freely generated $\infty$-groupoid. Unlike regular types, HITs
come with the ability to define equalities between terms. We will only introduce the syntax of HITs informally and refer the reader to \cite{lumsdaine2017semantics} for the detailed semantic interpretation (without which we would have no way of knowing our work was consistent). As an introductory
example, let us consider defining the circle, $S^{1}$, in homotopy type
theory. In this case the circle is generated by
\begin{enumerate}
\renewcommand{\theenumi}{(\arabic{enumi})}
\item a point $\base:S^{1}$,
\item a path $\loopme:\base=_{S^{1}}\base$.
\end{enumerate}
In other words the circle is a $\infty$-groupoid generated by a single
point $\base$ and a single path $\loopme$. Recall that in this context
equalities are thought of as paths. The question now becomes, what
is the induction principle of this type? The recursion principle is
not too hard to guess. For a type $C,$ to define a function $F:S^{1}\rightarrow C$
we must have the following:
\begin{enumerate}
\renewcommand{\theenumi}{(\arabic{enumi})}

\item $\base':C$,
\item $\loopme':\base'=_{C}\base'$.
\end{enumerate}
So, we must specify where the generators go to define a function from
$S^{1}$ to another $\infty$-groupoid. As we would expect $F(\base)\equiv\base'$
and $F(\loopme)=\loopme'$. Note the last equality is propositional;
this theory (as described in the HoTT book \cite{hottbook}) does
not compute, there is no canonical form for the path type (more importantly, the natural number type also has no canonical form, thought much progress has been made on this point \cite{shulman2015univalence, cohen2016cubical}). The induction principle is slightly more complicated.
Consider $C:S^{1}\rightarrow\text{Type}$, this can be thought of
as a fibration over $S^{1}$ (for each point in $s:S^{1}$ we have
a space $C(s))$, so a function $F:\Pi(t:S^{1}).C(t)$ is a section
of the fibration. Thus to define $F$ we require a point 

\[\base':C(\base),\]
and we'll need a \emph{dependent }path over $\loopme$
\[\loopme':\transport(\loopme)(\base')=\base.\]
Perhaps this is best understood by considering the total space of
the fibration, which here is just the product type $\Sigma(t:S^{1}).C(t)$.
A section $\Pi(t:S^{1}).C(t)$ is then really just a function $S^{1}\rightarrow\Sigma(t:S^{1}).C(t)$
\emph{where the first coordinate is given by the identity}. But we
can see how to define such a function using the recursion principle
and the definition of the product type $\Sigma$. We would need some
$\base':C(\base)$ so that we can get $(\base,\base'):\Sigma(t:S^{1}).C(t)$,
and we need $\loopme':\transport(\loopme)(\base')=\base$ so that
we can get $(\loopme,\loopme'):(\base,\base')=(\base,\base')$. 

We will mainly be dealing with pushouts (not surprising, since most
spaces are presented as gluing constructions which are in essence
iterated pushouts). A pushout can be thought of as two spaces glued
together in some way. In HoTT this gluing is accomplished by identifying
points with an equality. The classical version of what we are defining
is known as the double mapping cylinder. For spaces $X,Y,Z$ and continuous
maps $f\in Z\rightarrow X$ and $g\in Z\rightarrow Y$, we define
the double mapping cylinder as $Z\times I\sqcup X\sqcup Y$ quotiented
by the relation $(z,0)\simeq f(z)$ and $(z,1)\simeq g(z)$, so for
each point in $z$ there is a path between $f(z)$ and $g(z)$. Fix
types $X,Y,Z$. Given functions $f:Z\rightarrow X$ and $g:Z\rightarrow Y$,
we define the pushout $X+_{Z}Y$ by the generators
\begin{enumerate}
\renewcommand{\theenumi}{(\arabic{enumi})}
\item $\inl:X\rightarrow X+_{Z}Y$,
\item $\inr:Y\rightarrow X+_{Z}Y$ and,
\item $\glue:\Pi(z:Z).\inl(f(z))=\inr(g(z))$.
\end{enumerate}
Intuitively this says we have a copy of $X$ and a copy of $Y$ and
for every $z:Z$ we set $f(z)$ equal to $g(z)$. The recursion principle
is as follows. For a type $C$, to define a function $F:X+_{Z}Y\rightarrow C$
we must of course have functions 
\[\inl':X\rightarrow C
\text{and}\inr':Y\rightarrow C,\]
but since $\inl(f(z))=\inr(g(z))$ we also would expect to need an
equality 
\[\glue':\Pi(z:Z)\inl'(f(z))=\inr'(g(z)).\]
Of course $F(\inl(x))\equiv\inl'(x)$, $F(\inr(y))\equiv\inr'(x)$
and $F(\glue(z))=\glue'(z)$. The induction principle is defined similarly.
Given $C:X+_{Z}Y\rightarrow\text{Type}$, to define $F:\Pi(t:X+_{Z}Y).C(t)$
we need functions $\inl':\Pi(x:X).C(\inl(x))$, $\inr':\Pi(y:Y).C(\inr(y))$, 
together with a dependent path over $\glue$, namely
\[\glue':\Pi(z:Z).\transport(\glue)(\inl'(f(z))=\inr'(g(y)).\]
Note that in homotopy type theory we do not need to worry about defining
topology or any such thing. Indeed our definition looks like the definition
of set pushout, and yet we are really dealing with $\infty$-groupoids.

\subsection{Preliminaries}

Unless otherwise stated, all our definitions will be for pointed types,
pointed functions, pointed isomorphisms, and so on. We will often
avoid explicitly showing that our definitions are pointed but in each
case it should be straightforward to do so. For a type pointed $X$
we will denote the point by $x_{0}$, similarly $y_{0}$ for $Y$,
etc.

We will also often not explicitly work out the group structure. For
example we may show two types are isomorphic and avoid explicitly
showing the maps to be homomorphisms. Again, it should be straightforward
for the reader to work out the details.

We use the name of isomorphisms, say $f:X\simeq Y$, to refer to the
function $X\rightarrow Y$ that belongs to the definition of $X\simeq Y$.

We will make use of many basic higher inductive types, all details
can be found in the HoTT book \cite{hottbook}. We will use four
different particular types of pushouts: suspension, wedge, smash product,
and cofibers. In the classical setting this would be like defining
suspension, wedge, smash product and the mapping cylinder using double
mapping cylinders.

We define the suspension $\Sigma X$ by the constructors $N:\Sigma X$,
$S:\Sigma X$ and $\merid:\Pi(x:X).N=S$. As we have seen this says
we have two elements $N$ and $S$ ('north and south pole') and for
every element of $X$ we have a path between $N$ and $S$. For example
the circle is the suspension of the two element type $S^{0}$, i.e.,
we have two points $N,S$ and two paths between them. Classically
this corresponds to the double mapping cylinder $X\times I\sqcup\{N\}\sqcup\{S\}$
quotiented by the relation $(x,0)\simeq N$ and $(x,1)\simeq S.$
We leave it to the reader to understand the remaining definitions
in terms of their classical analogues. 

We define the wedge $X\lor Y$ by the constructors $\leftme:X\rightarrow X\lor Y$,
$\rightme:Y\rightarrow X\lor Y$, and $\wglue:\leftme(x_{0})=\rightme(y_{0})$.

We define the smash product $X\land Y$ by the constructors $\smbase:X\land Y$,
$\smin:X\times Y\rightarrow X\land Y$ and $\smglue:\Pi(p:X\lor Y).\smbase=\smin(f(p))$
where $f$ is the canonical map $X\lor Y\rightarrow X\times Y$ defined
by mapping $\leftme(x)$ to $(x,y_{0})$, $\rightme(y)$ to $(x_{0},y)$
and $\wglue$ to $\refl_{(x_{0},y_{0})}$. We often write $\smglue(x)$
for $\smglue(\leftme(x))$ and $\smglue(y)$ for $\smglue(\rightme(y))$
when we feel it is clear from context.

Finally for a function $f:X\rightarrow Y$ we define the cofiber $C_{f}$
by constructors $\cfbase:C_{f}$, $\cfcod(f):Y\rightarrow C_{f}$
and $\cfglue(f):\Pi(x:X).\cfbase=\cfcod(f)(f(x))$. in the later two
cases we will sometimes omit the $f$ when it is clear from context. 

We also make use of the circle $S^{1}$ defined by constructors $\base:S^{1}$
and $\loopme:\base=\base$ and set truncation $\|X\|_{0}$ defined
by constructors $|-|_{0}:X\rightarrow\|X\|_{0}$ and $\texttt{isset}:\Pi(y:\|X\|_{0}).$ 

Recall the path space $\Omega X$ is defined by $x_{0}=x_{0}$.
We define $\Sigma^{k}X$ inductively as $\Sigma(\Sigma^{k-1}X)$,
whereas we define $\Omega^{k}X$ as $\Omega^{k-1}(\Omega X)$. This
is worth knowing as occasionally we make use of the definitional equality.
\begin{lemma}
(Functor Lemma) \label{lem:(Functor-Lemma)}Let $K$ be some fixed
type. The following can be given actions on functions and become functors
from pointed types to pointed types: $\Sigma -$, $\Omega-$, $\|-\|_{0}$,$-\land K$,
$K\rightarrow -$. Moreover
 $\pi_{k} (-)
 \coloneqq
 \|\Omega^{k} -\|_{0}$
is a functor from pointed types to pointed sets, to groups for $k\geq1$
and to abelian groups for $k\geq2$.
\end{lemma}
\begin{proof}
This is fairly straightforward and much of it as been proven elsewhere,
see \cite{cavallo2015synthetic,hottbook}.
\end{proof}
\begin{definition}
\label{def:homologytheory}A collection of functors $H_{n}$ from
pointed types to groups is a homology theory when the following axioms
are satisfied.
\end{definition}
\begin{enumerate}
\renewcommand{\theenumi}{(\arabic{enumi})}

\item (Suspension) There exists an isomorphism $\susp:H_{n}(X)\simeq H_{n+1}(\Sigma X)$
and moreover this is natural, namely given $f:X\rightarrow Y$ we
have $\susp\circ(H_{n+1}(\Sigma f))\simeq H_{n}(f)\circ\susp$.
\item (Exactness) For any $f:X\rightarrow Y$ we have the following diagram is an exact sequence:

\end{enumerate}
\begin{center}
\begin{tikzcd} 
H_n (X) \arrow[r, "H_n(f)"] & H_n(Y)\arrow[rr, "H_n (\text{cfcod(f))}"] & &  H_n (C_f)
\end{tikzcd}
\end{center}

\section{Stable Homotopy Theory}

In this section we will define the stable homotopy groups in HoTT
and show that they satisfy the axioms of homology.
\begin{lemma}
\textbf{\label{lem: funmerid} }For every $X$ and $n$ there is a
function $\pi_{n}(X)\rightarrow\pi_{n+1}(\Sigma X)$.
\end{lemma}
\begin{proof}
This argument comes from chapter 8 of the HoTT book. We can define
$\merloop(X):X\rightarrow\Omega\Sigma X$ by sending $x$ to $\merid(x)\circ\merid(x_{0})^{-1}$,
then the required map is simply $\pi_{n}(\merloop(X))$.
\end{proof}
\begin{definition}
Given a sequence of types $X_{i}$ and maps $f_{i}:X_{i}\rightarrow X_{i+1}$
we define the colimit $\colim_{i\rightarrow\infty}X_{i}$ as a higher
inductive type with point constructors $\texttt{in}_{i}:X_{i}\rightarrow\colim_{i\rightarrow\infty}X_{i}$
and path constructors $\text{\texttt{coglue}}_{i}:\Pi(x:X_{i})\texttt{in}_{i}(x)=f_{i}\circ\text{\texttt{in}}_{i+1}(x)$.
\end{definition}
\begin{definition}
For a type $X$ the stable homotopy groups are defined as 
\[
\pi_{n}^{s}(X)\coloneqq\colim_{i\rightarrow\infty}\pi_{n+i}(\Sigma^{i}X),
\]
where the map from $\pi_{n+i}(\Sigma^{i}X)\rightarrow\pi_{n+i+1}(\Sigma^{i+1}X)$
is given by Lemma \ref{lem: funmerid}, namely \[\pi_{n+i}(\merloop(\Sigma^{i}X)).\]

\end{definition}
We denote $\pi_{n+i}(\merloop(\Sigma^{i}X))$ by $\phi_{i}(n,X)$,
where we may drop the $n$ and $X$ if it is clear from context. For
$j>i$ the composition $\phi_{i}\circ\phi_{i+1}...\circ\phi_{j-1}:\pi_{n+i}(\Sigma^{i}X)\rightarrow\pi_{n+j}(\Sigma^{j}X)$
will be denoted by $\Phi_{i}^{j}(n,X)$.

As the name suggests this sequence is eventually constant. To prove
this we will now recall some concepts and results from the homotopy
type theory book.
\begin{definition}
We say a function $f:X\rightarrow Y$ is $n-$connected provided
\[
\|\fib_{f}(y)\|_{n}
\]
is contractible for all $y:Y$. Recall $\fib_{f}(y)=\Sigma(x:X).f(x)=y$.
A type $X$ is $n-$connected provided $\|A\|_{n}$ is contractible.
\end{definition}
\begin{proposition}
\label{prop:(Freudenthal-Suspension-Theorem)}(Freudenthal Suspension
Theorem) Suppose $X$ is $n-connected$ then the map $\merloop(X)$
is $2n-$connected. \hfill \usebox{\proofbox}
\end{proposition}

\begin{proposition}
\label{prop:connpi}(8.8.5 from HoTT book) If $f:X\rightarrow Y$
is $n-connected$ and $k\leq n$ then $\pi_{k}(f)$ is an isomorphism. \hfill \usebox{\proofbox}
\end{proposition}

\begin{proposition}
\label{prop:suspconn}(8.2.1 from HoTT book) If $X$ is $n-$connected
then $\Sigma X$ is $n+1$-connected. \hfill \usebox{\proofbox}
\end{proposition}
\begin{lemma}
Suppose we have a sequence of types $X_{i}$ and functions $f_{i}:X_{i}\rightarrow X_{i+1}$
and suppose there exists an $i_{0}$ such that for all $i\geq i_{0}$,
$f_{i}$ is an isomorphism, then 
\[
\colim_{i}X_{i}\simeq X_{i_{0}}.
\]
\end{lemma}
\begin{proof}
We'll start by constructing a map $\colim_{i}X_{i}\rightarrow X_{i_{0}}$.
For each $x:X_{i}$ we want an element of $X_{i_{0}}$. If $i<i_{0}$
we use $f_{i_{0}-1}(f_{i_{0}-2}(...f_{i}(x))...))$. If $i>i_{0}$
we do a similar trick with $f_{i}^{-1}$'s. Finally if $i=i_{0}$
we merely use $x$. It then becomes trivial to show that our maps
satisfy the required equalities and hence we get a map $\colim_{i}X_{i}\rightarrow X_{i_{0}}$
as desired. On the other hand we use the canonical map $X_{i_{0}}\rightarrow\colim_{i}X_{i}$.
It is not difficult to see that these two maps form an inverse. 
\end{proof}
\begin{theorem}
(Stability) Given a type $X$ and $n:\mathbb{N}$ there exists an
$i_{X,n}$ such that each $\phi_{i}(n,X)$ is an isomorphism for $i\geq i_{X,n}$.
Moreover this will imply that there is an isomorpshim 
\[
\stab(n,X):\pi_{n}^{s}(X)\simeq\pi_{n+i_{X,n}}(\Sigma^{i_{X,n}}X).
\]
\end{theorem}
\begin{proof}
Proposition \ref{prop:suspconn} tells us that if $X$ is $n-$connected
then $\Sigma X$ is $n+1$-connected. Therefore as soon as $i$ is
large enough so that $2i\geq n+i$ then $\merloop(\Sigma^{i}X)$ will
be $2i$ connected by Prop \ref{prop:(Freudenthal-Suspension-Theorem)}
and hence $\pi_{n+i}(\merloop(\Sigma^{i}X))$ will be an isomorphism
by Prop \ref{prop:connpi}. The second result follows from the above
Lemma.
\end{proof}
The stability theorem will let us avoid working explicitly with colimits
when proving facts about the stable homotopy groups. Indeed one could
take the above theorem as the definition of the stable homotopy groups.
\begin{corollary}
Consider $X,\ n,\ F,\ i_{X,n}$ defined as in the previous result.
For any $i\geq i_{X,n}$ we can show $\pi_{n}^{s}(X)\simeq\pi_{n+i}(\Sigma^{i}X).$
\end{corollary}
\begin{proof}
This is straightforward, first invoke the previous Theorem to get
$\stab(n,X):\pi_{n}^{s}(X)\simeq\pi_{n+i_{X}}(\Sigma^{i_{X}}X)$ we
can then get $\pi_{n+i_{X}}(\Sigma^{i_{X}}X)\simeq\pi_{n+i}(\Sigma^{i}X)$
by composing with $\Phi_{i_{X}}^{i}$.
\end{proof}
We wish to show that $\pi_{n}^{s}(-)$ is a homology theory but first
we must show it is a functor. To that end we start by defining its
action on functions.
\begin{definition}
Given an $f:X\rightarrow Y$ we let $k=\max(i_{X},i_{Y})$, Note that
from Lemma \ref{lem:(Functor-Lemma)} we already have $\pi_{n+k}(\Sigma^{k}-)$
is a functor.\textbf{ }We define $\pi_{n}^{s}(f)$ as 
\[
(\stab(n,X)\circ\Phi_{i_{x}}^{k})\circ\pi_{n+k}(\Sigma^{k}f)\circ(\stab(n,X)\circ\Phi_{i_{y}}^{k})^{-1}.
\]
\end{definition}
Note that it is trivial to show $\pi_{n}^{s}(id)=id$ but showing
that $\pi_{n}^{s}$ preserves composition is not quite so easy. We
require the following lemma which allows us to replace $k$ in the
definition of $\pi_{n}^{s}(f)$.
\begin{lemma}
\label{lem:insertion}Let $k$ be as above then for any $k'\geq k$
we have 
\[
\pi_{n}^{s}(f)=(\stab\circ\Phi_{i_{x}}^{k'})\circ\pi_{n+k'}(\Sigma^{k'}f)\circ(\stab\circ\Phi_{i_{y}}^{k'})^{-1}.
\]
\end{lemma}
\begin{proof}
Clearly by definition and simple properties of repeated composition
it suffices to show
\[
\pi_{n+k}(\Sigma^{k}f)=\Phi_{k}^{k'}\circ\pi_{n+k'}(\Sigma^{k'}f)\circ(\Phi_{k}^{k'})^{-1}.
\]
Its also easy to see that it suffices to prove this for $k'=k+1$.
So we need the following square to commute:

\begin{center}
\begin{tikzcd} 
\pi_{n+k}(\Sigma^k X) \arrow[rrr, "\pi_{n+k}(\Sigma^k f)"] \arrow[d, "\phi_k(X)" ] 
& & &  \pi_{n+k}(\Sigma^k Y)  \arrow[d, "\phi_k(Y)" ] \\ 
\pi_{n+k+1}(\Sigma^{k+1} X) \arrow[rrr, "\pi_{n+k+1}(\Sigma^{k+1} f)" ] 
& & & \pi_{n+k+1}(\Sigma^{k+1} Y) 
\end{tikzcd}
\end{center}
But by the functoriality Lemma above we can reduce this to showing
the following commutes:

\begin{center}
\begin{tikzcd}  
\Sigma^k X \arrow[rr, "\Sigma^k f"] \arrow[d ] 
& & \Sigma^k Y \arrow[d ] \\ 
\Omega \Sigma^{k+1}X \arrow[rr, "\Omega\Sigma^{k+1} f" ] 
& & \Omega \Sigma^{k+1} Y
\end{tikzcd}
\end{center}
where the left and right arrows are $\merloop(\Sigma^{k}X)$ and $\merloop(\Sigma^{k}Y)$
respectively. In fact it will be clearer to prove a more general statement
namely for any $f:X\rightarrow Y$ the following commutes:
\begin{center}
\begin{tikzcd} 
X \arrow[r, "f"] \arrow[d ] 
& Y \arrow[d ] \\ 
\Omega \Sigma X \arrow[r, "\Omega\Sigma f" ] 
& \Omega \Sigma Y
\end{tikzcd}
\end{center}
This is easy to show. For any $x:X$ we wish to show 
\[
\merid(f(x))\circ\merid(y_{0})^{-1}=\apme_{\Sigma f}(\merid(x)\circ\merid(x_{0})^{-1})
\]
but the right hand side is simply 
\[
\merid(f(x))\circ\merid(f(x_{0}))^{-1},
\]
and of course $f(x_{0})=y_{0}$, completing the proof.
\end{proof}
\begin{theorem}
$\pi_{n}^{s}(-)$ is a functor. In particular it remains to show given
$f:X\rightarrow Y$ and $g:Y\rightarrow Z$ we have $\pi_{n}^{s}(g\circ f)=\pi_{n}^{s}(g)\circ\pi_{n}^{s}(f)$.
\end{theorem}
\begin{proof}
Let $k_{1}=\max(i_{x},i_{y})$ $k_{2}=\max(i_{y},i_{z})$ and $k_{3}=\max(i_{x},i_{z}).$
Written out in full our goal is to show that
\[
\stab(X)\circ\Phi_{i_{x}}^{k_{1}}(X)\circ\pi_{n+k_{1}}(\Sigma^{k_{1}}f)\circ(\stab(Y)\circ\Phi_{i_{Y}}^{k_{1}}(Y))^{-1}
\]
 composed with 
\[
\stab(Y)\circ\Phi_{i_{x}}^{k_{2}}(Y)\circ\pi_{n+k_{2}}(\Sigma^{k_{2}}g)\circ(\stab(Z)\circ\Phi_{i_{Z}}^{k_{2}}(Z))^{-1}
\]
is equal to
\[
\stab(X)\circ\Phi_{i_{x}}^{k_{3}}(X)\circ\pi_{n+k_{3}}(\Sigma^{k_{3}}(f\circ g))\circ(\stab(Z)\circ\Phi_{i_{Z}}^{k_{3}}(Z))^{-1}.
\]
But applying the above Lemma to each of these components we can replace
$k_{1},k_{2},k_{3}$ with $k\coloneqq\max(k_{1},k_{2},k_{3})$. Once
this is done the equation can be simplified until it suffices to show
\[
\pi_{n+k}(\Sigma^{k}f)\circ\pi_{n+k}(\Sigma^{k}g)=\pi_{n+k}(\Sigma^{k}f\circ g),
\]
which follows from functoriality of $\pi_{n+k}$ and $\Sigma^{k}$
\end{proof}
\begin{theorem}
The functors $\pi_{n}^{s}$ satisfy the suspension axiom of homology
theories. 
\end{theorem}
\begin{proof}
We first construct $\susp(n,X):\pi_{n}^{s}(X)\simeq\pi_{n+1}^{s}(\Sigma X)$.
Let $k=\max(i_{X,n},i_{\Sigma X,n+1})$ then we have 
\[
\stab(n,X)\circ\Phi_{i_{x}}^{k+1}:\pi_{n}^{s}(X)\simeq\pi_{n+(k+1)}(\Sigma^{k+1}X).
\]
Similarly we have 
\[
\stab(\Sigma X,n+1)\circ\Phi_{i_{\Sigma X}}^{k}:\pi_{n+1}^{s}(\Sigma X)\simeq\pi_{(n+1)+k}(\Sigma^{k}\Sigma X).
\]
Thus it suffices to find 
\[
\susp_{0}(X,n,k):\pi_{n+(k+1)}(\Sigma^{k+1}X)\simeq\pi_{(n+1)+k}(\Sigma^{k}\Sigma X).
\]
We define $\susp_{0}$ by first transporting along the equality $n+(k+1)=(n+1)+k$
and then transporting along the equality$\Sigma^{k+1}X=\Sigma^{k}\Sigma X$.
Therefore (omitting many dependencies) we have 
\[
\susp(n,X)\coloneqq\stab(X)\circ\Phi_{i_{x}}^{k+1}\circ\susp_{0}\circ(\Phi_{i_{\Sigma X}}^{k})^{-1}\circ\stab(\Sigma X)^{-1}.
\]
We will need to show that we can replace $k$ with any $k'>k$ as
we did for $\pi_{n}^{s}(f)$. In this case, as before, begin by performing
the obvious cancellations and reduce to the case $k'=k+1$ then it
will suffice to show the following commutes

\begin{center}
\begin{tikzcd}  
\pi_{n+(k+1)} (\Sigma^{k+1} X) \arrow[rrr, "\susp_0 (X \, n \,k)"] \arrow[d,"\phi_{k+1} (n \, X)"] 
& & & \pi_{(n+1)+k}(\Sigma^k \Sigma X) \arrow[d, "\phi_{k} (n+1 \, \Sigma X)" ] \\ 
\pi_{(n+(k+1))+1} \Sigma^{(k+1)+1} X \arrow[rrr, "\susp_0(X \, n+1 \, k+1)"] 
& & & \pi_{((n+1)+k)+1}(\Sigma^{k+1} \Sigma X)
\end{tikzcd}
\end{center}
Note, since $n+(k+1)\equiv(n+k)+1$ and $\Sigma^{k+1}X\equiv\Sigma\Sigma^{k}X$,
we see that both the top and bottom are merely transporting along
the equalities $(n+(k+1))=(n+1)+k$ and $\Sigma^{k+1}X=\Sigma^{k}\Sigma X$.
Moreover $\phi_{k+1}(n,X)$ is $\pi_{n+(k+1)}(\merloop(\Sigma^{k+1}X))$
and $\phi_{k}(n+1,\Sigma X)$ is $\pi_{(n+1)+k}(\merloop(\Sigma^{k}\Sigma X))$,
notice that they depend on the same terms as the types in the corners.
Therefore filling this square becomes a simple exercise in using path
induction.

Now that we know we can increase $k$ in the definition above, we
can show the $\susp$ is natural. We must show
\[
\susp(n,X)\circ\pi_{n+1}^{s}(\Sigma f)=\pi_{n}^{s}(f)\circ\susp(n,Y).
\]
By expanding we see that we must show
\begin{align*}
\stab(n,X)\circ\Phi_{i_{x}}^{k_{1}+1}\circ\susp_{0}(X,n,k_{1})\circ(\Phi_{i_{\Sigma X}}^{k_{1}})^{-1}\circ\stab(n+1,\Sigma X)^{-1}\\
\circ(\stab(n+1,\Sigma X)\circ\Phi_{i_{\Sigma X}}^{k_{2}})\circ\pi_{(n+1)+k_{2}}(\Sigma^{k_{2}}\Sigma f)\circ(\stab(n+1,\Sigma Y)\circ\Phi_{i_{\Sigma Y}}^{k_{2}})^{-1}
\end{align*}
is equal to
\begin{align*}
(\stab(n,X)\circ\Phi_{i_{x}}^{k_{3}})\circ\pi_{n+k_{3}}(\Sigma^{k_{3}}f)\circ(\stab(n,Y)\circ\Phi_{i_{Y}}^{k_{3}})^{-1}\circ\\
\stab(n,Y)\circ\Phi_{i_{Y}}^{k_{4}+1}\circ\susp_{0}(Y,n,k_{4})\circ(\Phi_{i_{\Sigma Y}}^{k_{4}})^{-1}\circ\stab(n+1,\Sigma Y)^{-1}.
\end{align*}
Replace everything so that $k_{3}=k_{4}+1=k_{2}+1=k_{1}+1$, we are
then left with showing 
\[
\susp_{0}(X,n,k)\circ\pi_{(n+1)+k}(\Sigma^{k}\Sigma f)=\pi_{n+k+1}(\Sigma^{k+1}f)\circ\susp_{0}(Y,n,k).
\]
Split this into two squares. The first we transport along the equality
$(n+1)+k=n+k+1$ and so easily get that it commutes. The second we
transport along the equality $\Sigma^{k+1}X=\Sigma^{k}\Sigma X$.
After applying the functoriality of $\pi_{n+k+1}$ it suffices to
show the following commutes:

\begin{center}
\begin{tikzcd}  
\Sigma^{k+1}X \arrow[r] \arrow[d,"\Sigma^{k+1}f"] 
& \Sigma^k \Sigma X \arrow[d, "\Sigma^k \Sigma f" ] \\ 
\Sigma^{k+1}Y \arrow[r] 
& \Sigma^k \Sigma Y 
\end{tikzcd}
\end{center}
We do this by induction. For $k=1$ this is trivial, otherwise if
it is true for $k-1$ recall that we get the equality $\Sigma^{k+1}X=\Sigma^{k}\Sigma X$
by rewriting $\Sigma^{k+1}X=\Sigma\Sigma^{k}X$ and $\Sigma^{k}\Sigma X=\Sigma\Sigma^{k-1}\Sigma X$
and then by using the fact that $k=(k-1)+1$ and applying $\Sigma$
to the equality for $k-1$. Thus the entire square can be again be
decomposed. First we transport along $k=(k-1)+1$ then we have simply
$\Sigma$ applied to the $k-1$ square. Both these new squares are
easy to fill.
\end{proof}
It remains to show that $\pi_{n}^{s}(-)$ satisfies the exactness
axiom. To this end we need the following facts.
\begin{theorem}
\label{thm:Black Massey}(Blakers-Massey Theorem) if $A$ is $n-$connected
and $f:A\rightarrow B$ is $m$ connected then $\pi_{\ell}(A)\rightarrow\pi_{\ell}(B)\rightarrow\pi_{\ell}(C_{f})$
is exact for $\ell \leq n+m$.
\end{theorem}
\begin{proof}
This is a straightforward consequence of the version of Blakers-Massey
Theorem proven in \cite{finster2016mechanization}.
\end{proof}

We will need to know that suspension preserves cofibers. In fact later
we will need to know that smash product preserves cofibers. So we
first prove that suspension is a special case of smash and then we
show the more general result. Classically this follows since smash is left adjoint to pointed arrow and hence preserves colimits. However this general statement is not easy to prove in HoTT and so instead we prove the result directly.
\begin{lemma}
\label{lem:suspsmash}there exists isomorphisms $\Sigma X\simeq S^{1}\land X$,
$\suspsm(X,K):\Sigma(X\land K)\simeq(\Sigma X)\land K$, $\suspsmtwo(X,K):\Sigma(X\land K)\simeq X\land(\Sigma K)$,
moreover this is natural, in particular we will use the following
facts for $f:X\rightarrow Y,\ g:K\rightarrow L$
\[
\suspsm(X,K)\circ(\Sigma f)\land K=\Sigma(f\land K)\circ\suspsm(Y,K),
\]

\[
\suspsmtwo(X,K)\circ(f\land\Sigma K)=\Sigma(f\land K)\circ\suspsm2(Y,K),
\]
\[
X\land g\circ f\land L=f\land K\circ Y\land g,
\]
\begin{align*}
\suspsmtwo(\Sigma X,K)\circ\Sigma X\land g\circ\suspsm(X,L)^{-1}\\
=\Sigma\big(\suspsm(X,K)^{-1}\circ\suspsmtwo(X,K)\circ(X\land g)\big).
\end{align*}
\end{lemma}
\begin{proof}
In \cite{brunerie2016homotopy} we see that $\Sigma X\simeq S^{1}\land X$.
Moreover Floris van Doorn has manged to show the smash product is a symmetric monoidal product (in
particular associative and communative in a natural way) \cite{van2018formalization}. Assuming the result we can prove the first isomorphism $\suspsm$ as follows
\[
\Sigma(X\land K)\simeq S^{1}\land(X\land K)\simeq(S^{1}\land X)\land K\simeq\Sigma X\land K.
\]
To prove the first equation we can use the symmetric monoidal structure
of $\land$, the only thing missing is naturality of $\Sigma X\simeq S^{1}\land X$
but this can be easily deduced as follows. Let $\Phi(X):\Sigma X\rightarrow S^{1}\land X$
be the isomorphism given in \cite{brunerie2016homotopy}. Let $f:X\rightarrow Y$,
we wish to show the following commutes

\begin{center}
\begin{tikzcd}  
\Sigma X \arrow[r, "\Sigma f"] \arrow[d,"\Phi (X)"] 
& \Sigma Y \arrow[d, "\Phi (Y)" ] \\ 
S^{1} \land X \arrow[r, "S^{1} \land f"] 
& S^{1} \land  Y 
\end{tikzcd}
\end{center}
That is to say we want to show for $q:\Sigma X$ that 
\[
\Sigma f\circ\Phi(Y)(q)=\Phi(X)\circ S^{1}\land f(q).
\]
We proceed by induction on $q:\Sigma X$ for $q\equiv N$ we trace
through the definitions and see we just need to show $\cfbase=\cfbase$
which we can do by $\refl_{\cfbase}$. Similarly for $q\equiv S$.
Finally for $q=\merid(x)$ it will suffices to show 
\[
\Sigma f\circ\Phi(Y)(\merid(x))=\Phi(X)\circ S^{1}\land f(\merid(x)).
\]
Again we trace through the definitions and see it suffices to show
\begin{align*}
S^{1}\land f\big(\smglue(\leftme(x))\circ\apme_{\lambda s.\smin(s,x)}(\loopme)\circ\smglue(\leftme(x))^{-1}\big)\\
=\smglue(\leftme(f(x)))\circ\apme_{\lambda s.\smin(s,fx)}(\loopme)\circ\smglue(\leftme(f(x))^{-1},
\end{align*}
which is clear. The second isomorphism $\suspsmtwo$ and the second
equality can be shown in a similar manner. The third equality follows
directly from the symmetric monoid structure, and finally the last
equality looks more formidable but can also be easily deduced from
the above square and the properties in \cite{brunerie2016homotopy}.
\end{proof}
\begin{lemma}
\label{lem:SmRec}To construct a term $F:\Pi(x:X\land Y).P(x)$ it
suffices to give the following ingredients: 
\end{lemma}
\begin{enumerate}
\renewcommand{\theenumi}{(\arabic{enumi})}
\item $\smbase':P(\smbase).$
\item $\smin':\Pi(p:X\times Y).P(\smin(p)).$
\item $\smglueleft':\Pi(x:X)\smbase'=_{\smglue(\leftme x))}\smin'(x,y_{0}).$
\item $\smglueright':\Pi(y:Y)\smbase'=_{\smglue(\rightme(y))}\smin'(x_{0},y).$
\end{enumerate}
Moreover $F$ will satisfy $F(\smbase)\equiv\smbase'$, $F(\smin(p))\equiv\smin'(p)$,
\\$F(\smglue(\leftme(x))=\smglueleft'$ and $F(\smglue(\rightme(y)))=\smglueright'$.

Now let $f:X\rightarrow Y$ be fixed. To construct a term $F:\Pi(x:C_{f\land K}).P(x)$
it suffices to give the following ingredients:
\begin{enumerate}
\renewcommand{\theenumi}{(\arabic{enumi})}
\item $\cfbase':P(\cfbase)$.
\item $\cfcod':\Pi(p:Y\land K).P(\cfcod(p))$.
\item $\cfgluesmbase':\cfbase'=_{\cfglue(\smbase)}\cfcod'(f\land K(\smbase))$.
\item $\cfgluesmin':\Pi(q:X\times K)\smbase'=_{\cfglue(\smin(q))}\cfcod'(f\land K(\smin(q))$.
\end{enumerate}
Moreover $F$ will satisfy $F(\cfbase)\equiv\cfbase'$, $F(\cfcod(p))\equiv\cfcod'(p)$, \\$F(\cfglue(\smbase))=\cfgluesmbase'$ and $F(\cfglue(\smin(q)))=\cfgluesmin'$.
\begin{proof}
The first part follows from the fact that the smash product $X\land Y$ can be defined as the pushout of $2 \leftarrow X + Y \rightarrow X \times Y$. This is shown in the PhD of Floris van Doorn \cite{van2018formalization}. Using this definition it is clear that the given ingredients suffice to construct an element of $\Pi(x:X\land Y).P(x)$.
Similarly the second part will follow from the following claim: $C_{f\land K}$  can be defined as the pushout of $2 \leftarrow 1 + X \times K \rightarrow Y\land K$ where the map $1 + X\times K\rightarrow Y\land K$ is defined by sending the point to the basepoint of $Y\land K$ and sending $(x,k)$ to $f(\smin(x,k))$. To prove this claim consider the following diagram.
\begin{center}
\begin{tikzcd}  
1+X+K \arrow[r] \arrow[d] 
&  1 + X\times K \arrow[r] \arrow[d] 
& Y\land K \arrow[d] \\ 
1 + 2 \arrow[r] 
& 1 + X\land K \arrow[r] \arrow[d]
& Y \land K \arrow[d]  \\
& 2 \arrow[r]
& C_{f\land K}
\end{tikzcd}
\end{center}
We want to show te right rectangle is a pushout. It is not hard to show that the top left square, the bottom right square and the top rectangle are all pushouts. It hence follows from the pushout lemma that the top right square and hence the right rectangle are also pushouts.
\end{proof}
\begin{theorem}
There exists an isomorphism
\[
\smcf:C_{f\land K}\simeq C_{f}\land K,
\]
and $\cfcod(f\land K)\circ\smcf=\cfcod(f)\land K$.
\end{theorem}
\begin{proof}
We will make use of \ref{lem:SmRec} whenever we can. 

Define $\phi:C_{f\land K}\rightarrow C_{f}\land K$ by sending $\cfbase_{f\land K}$
to $\smbase$, $\cfcod_{f\land K}(q)$ to $\cfcod\land K(q)$ then
for $p:A\land K$ we need to show 
\[
\smbase=\cfcod_{f}\land K(f\land K(p)).
\]
For $p\equiv\smbase$ we can just use $\refl_{\smbase}$. For $p\equiv\smin(a,k)$
we need to show 
\[
\smbase=\smin(\cfcod(f(a)),k),
\]
which can be done by $\smglue(k)\circ\smin(\cfglue(a),k)$.\textbf{ }

Next we construct a map $\psi:C_{f}\land K\rightarrow C_{f\land K}.$
Map $\smbase$ to $\cfbase$. For $\smin(t,k)$ with $t:C_{f}$ we
induct on $t.$ For $t\equiv\cfbase_{f}$ use $\cfbase_{f\land K}$,
for $t\equiv\cfcod_{f}(b)$ we use $\cfcod_{f\land K}(\smin(b,k))$
and finally for $a:A$ we must show 
\[
\cfbase_{f\land K}=\cfcod(\smin(f(a),k)).
\]
But since $\smin(f(a),k)\equiv f\land K(\smin(a,k))$ we can simply
apply $\cfglue_{f\land K}(a,k)$.\textbf{ }For $k:K$ we need that
$\cfbase_{f\land K}=\cfbase_{f\land K}$ which we of course do with
$\refl$, and for $t:\cfbase_{f}$ we must show
\[
\cfbase=\psi(\smin(t,k_{0})),
\]
 ($\psi$ standing in for the partially defined function above), which
we prove by induction on $t.$ For $t\equiv\cfbase$ we need $\cfbase=\cfbase$
which we do by $\refl$. For $t\equiv\cfcod_{f}(b)$ we need to show
\[
\cfbase=\cfcod(\smin(b,k_{0})).
\]
Now we have $\smglue(\leftme(b)):\smin(b,k_{0})=\smbase$ and $\cfglue(\smbase):\cfbase=\cfcod(\smbase)$
so we simply use $\cfglue(\smbase)\circ\apme_{\cfcod}(\smglue(\leftme(b)))$.
For the case of $t=\cfglue(a)$ this reduces to showing 
\[
\cfglue(\smbase)\circ\apme_{\cfcod}(\smglue(\leftme(f(a)))\circ\psi(\glue(a))^{-1}=\refl_{\cfbase},
\]
which is 
\[
\cfglue(\smbase)\circ\apme_{\cfcod}(\smglue(\leftme(f(a))\circ\cfglue(\smin(f(a),k_{0}))^{-1}=\refl_{\cfbase}.
\]
But $\smglue(\leftme(f(a))=\apme_{f\land K}(\smglue(\leftme(a))$,
so we replace the middle with \[\apme_{f\land K\circ\cfcod}(\smglue(\leftme(a))\].
Note $f\land K\circ\cfcod$ is constant and the proof of this is simply
$\cfglue$, therefore by a simple path induction argument the above
equality is satisfied.

Let's now show $\Pi(z:C_{f\land K}).\psi(\phi(z))=z$. For $z\equiv\cfbase$
we use $\refl_{\cfbase}$ for $z\equiv\cfcod(y)$ for $y:B\land K$
we must induct on $y$. For $y\equiv\smbase$ we use the proof that
\[
\cfbase=\cfcod(\smbase),
\]
namely $\cfglue(\smbase)$. For $y\equiv\smin(b,k)$ we can just use
$\refl$. Now for $b:B$ we need to check
\[
\psi(\cfcod\land K(\smglue(b))^{-1}\circ\cfglue(\smbase)\circ\cfcod(\smglue(b))=\refl,
\]
but the first term is by definition equal to
\[
(\cfglue(\smbase)\circ\cfcod(\smglue(\leftme(b)))^{-1}.
\]
so this is clear. Next for $k:K$ we must show 
\[
\psi(\cfcod\land K(\smglue(k))^{-1}\circ\cfglue(\smbase)\circ\cfcod(\smglue(k))=\refl.
\]
but here the first term computes to 
\[
\psi(\smglue_{C_{f}\land K}(k)\circ\apme_{\smin(-,k)}((\cfglue(a_{0})\circ\apme_{\cfcod}(p)))^{-1}).
\]
where $p$ is the proof that $f(a_{0})=b_{0}$. which is 
\[
(\cfglue(\smin(a_{0},k))\circ\apme_{\lambda x.\cfcod(\smin(x,k))}(p))^{-1}.
\]
The proof can then be completed by noting that since $\cfcod(f\land K(x))$
is constant by $\cfglue$ we have 
\[
\cfcod(f\land K(\smglue(k))=\cfglue(\smbase)^{-1}\circ\cfglue(\smin(a_{0},k)),
\]
while on the other hand, by definition of $f\land K$,
\[
\cfcod(f\land K(\smglue(k))=\cfcod(\smglue(k))\circ\cfcod(\smin(p,k)),
\]
so we solve for $\cfcod(\smglue(k))$, plug it into our original equation
and everything easily cancels.

We have two more cases to check. First corresponding to $\cfglue(\smbase_{A\land K})$,
we would like the following
\[
\psi(\phi(\cfglue(\smbase))=\refl,
\]
however this is immediate from the definition. Next corresponding
to $\cfglue(\smin(a,k))$ we would like to show 
\[
\psi(\phi(\cfglue(\smin(a,k))=\cfglue(\smin(a,k)),
\]
but again this follows from the definition since the left hand side
is 
\[
\psi(\smglue(k)\circ\smin(\cfglue(a)),k),
\]
which is 
\[
\refl\circ\cfglue(\smin(a,k)).
\]
Finally we need to prove $\Pi(z:C_{f}\land K).\phi(\psi(z))=z$. We
first consider $z\equiv\smbase$ then we can just use $\refl$. If
$z\equiv\smin(t,k)$ we induct on $t.$ For $t\equiv\cfbase$ we must
show 
\[
\smbase=\smin(\cfbase,k),
\]
which we prove by $\smglue_{C_{f}\land K}(k)$. Then for $t\equiv\cfcod(b)$
we must show 
\[
\smin(\cfcod(b),k)\equiv\smin(\cfcod(b),k),
\]
which is done by $\refl.$ Finally we need to consider the case for
$t=\glue(a)$. In other words we must show 
\[
\phi(\psi(\smin(\glue(a),k)))\circ\smin(\cfglue(a),k)^{-1}\circ\smglue(k)^{-1}=\refl.
\]
The first term computes to 
\[
\phi(\cfglue(\smin(a,k)))=\smglue(k)\circ\smin(\cfglue(a),k),
\]
so the required equality is trivial.

Next consider the case $z=\smglue(k)$ where $k:K$. We must show
\[
\phi(\psi(\smglue(k)))\circ\smglue(k)\circ\smglue(k)^{-1}=\refl,
\]
so it suffices to show 
\[
\phi(\psi(\smglue(k))=\refl,
\]
which is by definition of $\psi$.

Next consider the case $z=\smglue(t)$, where $t:C_{f}$. We induct
on $t$. For $t\equiv\cfbase$ show 
\[
\phi(\psi(\smglue(\cfbase))\circ\smglue(k_{0})\circ\smglue(\cfbase)^{-1}=\refl.
\]
However $\smglue(k_{0})=\smglue(\cfbase)$, so we just need 
\[
\phi(\psi(\smglue(\cfbase))=\refl,
\]
which follows by definition. For $t\equiv\cfcod(b)$ we must show
\[
\phi(\psi(\smglue(\cfcod(b))=\smglue(\cfcod(b)),
\]
this follows from the definition since the left hand side is 
\[
\phi(\cfglue(\smbase)\circ\apme_{\cfcod}(\smglue(\leftme(b)))),
\]
which is 
\[
\refl\circ\cfcod\land K(\smglue(b))=\smglue(\cfcod(b)).
\]
Finally for $t=\glue(a)$ we must show some 3 dimensional condition.
Using terminology from \cite{licata2015cubical} we have already
filled the following squares

\begin{center}
\begin{tikzcd}  
\smglue \arrow[rrr, "\smglue (\cfbase)"] \arrow[d,"\refl"] 
& & & \smin(\cfbase, k_0) \arrow[d, "\smglue (k_0 )" ] \\ 
\phi (\psi (\smbase)) \arrow[rrr, "\phi (\psi (\smglue (\cfbase)))"] 
& & &\phi (\psi (\smin(\cfbase ,k_0)
\end{tikzcd}
\end{center}

\begin{center}
\begin{tikzcd}  
\smglue \arrow[rrr, "\smglue (\cfcod(b))"] \arrow[d,"\refl"] 
& & &\smin(\cfcod(b), k_0) \arrow[d, "\refl" ] \\ 
\phi (\psi (\smbase )) \arrow[rrr, "\phi (\psi (\smglue (\cfcod(b))))"] 
& & & \phi (\psi (\smin (\cfcod(b) ,k_0)
\end{tikzcd}
\end{center}
and we want a cube between them. We have $\cfglue(a):\cfbase=\cfcod(a).$
With some contemplation (and since square fillers are unique) it can
be seen that we are essentially in the following case: Suppose we
have types $X$, $A$, $f:X\rightarrow X$, $x:X$, $x':A\rightarrow X$,
$p:\Pi(a:A).x=x'(a)$, $\ell:f(x)=x$, $a_{1}:A$, $a_{2}:A$, and
finally $q:a_{1}=a_{2}$.

\begin{center}
\begin{tikzcd}  
x \arrow[r, "p(a_1)"] \arrow[d,"\ell"] 
&  x'(a_1) \arrow[d, "\text{filler}" ] \\ 
f(x) \arrow[r, "f(p(a_1))"] 
&  f(x'(a_1)
\end{tikzcd}
\end{center}

\begin{center}
\begin{tikzcd}  
x \arrow[r, "p(a_2)"] \arrow[d,"\ell"] 
&  x'(a_2) \arrow[d, "\text{filler}" ] \\ 
f(x) \arrow[r, "f(p(a_2))"] 
&  f(x'(a_2))
\end{tikzcd}
\end{center}
These two squares have a cube between them. This is clear since they
are the same square if path induction is applied to $q$.

It still remains to show the equality $\cfcod(f\land K)\circ\smcf=\cfcod(f\land K)$
but this follows quickly from the definition.
\end{proof}
\begin{corollary}
Given the above two results it follows that
\[
\suspcf:C_{\Sigma f}\simeq\Sigma C_{f},
\]
and $\cfcod(\Sigma f)\circ\suspcf=\Sigma\cfcod(f)$. \hfill\usebox{\proofbox}
\end{corollary}
\begin{corollary}
\label{cor:It-follows-that}It follows that $\Sigma^{k}$ preserves
cofibers in the same manner as $\Sigma$. \hfill \usebox{\proofbox}
\end{corollary}
\begin{lemma}
For any group homomorphism $f:X\rightarrow Y$ and $g:Y\rightarrow Z$,
if we have isomorphisms $\ell_{1}:X_{0}\rightarrow X$, $\ell_{2}:Y\rightarrow Y_{0}$
and $\ell_{3}:Z\rightarrow Z_{0}$ and then $\ell_{1}\circ f\circ\ell_{2}$,
$\ell_{2}^{-1}\circ g\circ\ell_{3}$ is exact if $f,g$ is exact.
\end{lemma}
\begin{proof}
Straightforward.
\end{proof}
\begin{theorem}
$\pi_{n}^{s}(-)$ satisfies the exactness axiom.
\end{theorem}
\begin{proof}
For any $f:X\rightarrow Y$ we wish to show that 

\begin{center}
\begin{tikzcd}  
\pi_n^s (X) \arrow[rr, "\pi_n^s f"] & & \pi_n^s (Y) \arrow[rr, "\pi_n^s \cfcod(f)"] & & \pi_n^s(C_f)
\end{tikzcd}
\end{center}
is exact. We start by expanding the definitions of $\pi_{n}^{s}(f)$
and $\pi_{n}^{s}(\cfcod(f))$ and using our replacement lemma as we've
done above. It then becomes immediately clear that we are in the same
position as the lemma above and so it suffices to show

\begin{center}
\begin{tikzcd}  
\pi_{n+k} (\Sigma^k X) \arrow[rr, "\pi_{n+k} (\Sigma^k f)"] & & \pi_{n+k} \Sigma^k(Y) \arrow[rrr, "\pi_{n+k} \Sigma^k \cfcod(f)"] & & & \pi_{n+k} \Sigma^k (C_f)
\end{tikzcd}
\end{center}
is exact where $k\geq i_{X},i_{Y},i_{C_{f}}$. But from Corollary
\ref{cor:It-follows-that} we can replace this with

\begin{center}
\begin{tikzcd}  
\pi_{n+k} (\Sigma^k X) \arrow[rr, "\pi_{n+k} (\Sigma^k f)"] & & \pi_{n+k} \Sigma^k(Y) \arrow[rrr, "\pi_{n+k} \cfcod(\Sigma^k f)"] & & & \pi_{n+k} (C_{\Sigma^k f}).
\end{tikzcd}
\end{center}
Now $\Sigma^{k}X$ is $k$ connected as long as $k\geq n$ (by \ref{prop:suspconn}).
Moreover $\Sigma^{k}Y$ is $k$ connected so the domain and codomain
of $\Sigma^{k}f$ is $k$ connected and hence $\Sigma^{k}f$ must
be $k-1$ connected (for a proof combine 7.5.6 and 7.5.11 from the
HoTT book). Therefore by the Blakers Massey Theorem (\ref{thm:Black Massey})
this sequence is exact.
\end{proof}

\begin{corollary}
\label{cor:stableishomology}$\pi_{n}^{s}(-)$ is a homology theory. \hfill \usebox{\proofbox}
\end{corollary}

\section{Extending the Result to Homology}

We first wish to extend the previous result by showing for a fixed
type $K$, $\pi_{n}^{s}(-\land K)$ is a homology theory.

It is clear that $\pi_{n}^{s}(-\land K)$ is a functor since it is
the composition of $\pi_{n}^{s}$ and $-\land K$ both of which we've
already shown are functors.
\begin{theorem}
$\pi_{n}^{s}(-\land K)$ satisfies the suspension axiom.
\end{theorem}
\begin{proof}
We start by constructing an isomorphism $\susptwo(n,X):\pi_{n}^{s}(X\wedge K)\simeq\pi_{n+1}^{s}((\Sigma X)\land K)$.
Of course we already have $\pi_{n}^{s}(X\land K)\simeq\pi_{n+1}^{s}(\Sigma(X\land K))$
by $\susp(n,X\land K)$. It therefore suffices to compose this with
the isomorphism $\pi_{n+1}^{s}(\suspsm):\pi_{n+1}^{s}(\Sigma(X\land K))\simeq\pi_{n+1}^{s}((\Sigma X)\land K))$.
Next we must show this is natural. We wish to show 
\[
\susptwo(n,X)\circ\pi_{n+1}^{s}((\Sigma f)\land K)=\pi_{n}^{s}(f\land K)\circ\susptwo(n,Y).
\]
Therefore, we must show
\begin{align*}
\susp(n,X\land K)\circ\pi_{n+1}^{s}(\suspsm(X))\circ\pi_{n+1}^{s}((\Sigma f)\land K)=\\
\pi_{n}^{s}(f\land K)\circ\susp(n,Y\land K)\circ\pi_{n+1}^{s}(\suspsm(Y).
\end{align*}
 Using naturality of $\susp$ we can simplify this to showing
\[
\pi_{n+1}^{s}(\suspsm(X))\circ\pi_{n+1}^{s}((\Sigma f)\land K)=\pi_{n+1}^{s}(\Sigma(f\land K)\circ\pi_{n+1}^{s}(\suspsm(Y)),
\]
which would follow provided we show 
\[
\suspsm(X)\circ(\Sigma f)\land K=\Sigma(f\land K)\circ\suspsm(Y).
\]
This is true by Lemma \ref{lem:suspsmash} above.
\end{proof}
\begin{theorem}
$\pi_{n}^{s}(-\land K)$ satisfies the exactness axiom. 
\end{theorem}
\begin{proof}
We wish to show

\begin{center}
\begin{tikzcd}  
\pi_n^s (X \land K) \arrow[rr, "\pi_n^s f\land K"] & & \pi_n^s (Y \land  K) \arrow[rr, "\pi_n^s \cfcod_f\land K"] & & \pi_n^s(C_f\land K)
\end{tikzcd}
\end{center}

is exact. Just like in the proof of exactness of $\pi_{n}^{s}$ we
can replace this with

\begin{center}
\begin{tikzcd}  
\pi_n^s (X \land K) \arrow[rr, "\pi_n^s f\land K"] & & \pi_n^s (Y \land  K) \arrow[rr, "\pi_n^s \cfcod_{f\land K}"] & & \pi_n^s(C_{f\land K}).
\end{tikzcd}
\end{center}

We are then done by exactness of $\pi_{n}^{s}$.
\end{proof}
\begin{corollary}
\label{cor:stablesm-is-homo}$\pi_{n}^{s}(-\land K)$ is a homology
theory. \hfill \usebox{\proofbox}
\end{corollary}

\subsection{Homology}

We are now ready to construct ordinary homology. 
\begin{definition}
A prespectrum is a sequence of types $K_{i}$ and maps $\kappa_{i}:\Sigma K_{i}\rightarrow K_{i+1}.$
\end{definition}
\begin{lemma}
Given a prespectrum $(K_{i},\kappa_{i})$ and a type $X$ we have a prespectrum
with types $X\land K_{i}$.
\end{lemma}
\begin{proof}
From Lemma \ref{lem:suspsmash} we have $\suspsmtwo(X,K_{i}):\Sigma(X\land K_{i})\rightarrow X\land\Sigma K_{i}$.
Compose this with $X\land\kappa_{i}$.
\end{proof}
\begin{lemma}
Given a prespectrum $(K_{i},\kappa_{i})$ we can construct a map \[\pi_{n+i}^{s}(K_{i})\rightarrow\pi_{n+i+1}^{s}(K_{i+1}).\] 
\end{lemma}
\begin{proof}
We first use $\susp(K_{i})$ from above to get $\pi_{n}^{s}(K_{i})\rightarrow\pi_{n+i+1}^{s}(\Sigma K_{i})$
then we compose with $\pi_{n+i+1}^{s}(\kappa_{i})$.
\end{proof}
\begin{theorem}
Using the functions we just computed above we define \[H_{n}(X)=\|\colim_{i}\pi_{n+i}^{s}(X\land K_{i})\|_{0}.\]
This is a homology theory. 
\end{theorem}
\begin{remark}
Truncated colimits of sets behave much like we expect form the classical
setting. In particular note the following. To define a function from
$\|\colim_{i}(X_{i},\phi_{i})\|_{0}\rightarrow\|\colim_{i}(Y,\psi_{i})\|_{0}$
where each $X_{i},Y_{i}$ is a set it suffices to find maps $f_{i}:X_{i}\rightarrow Y_{i}$
such that $f_{i}\circ\psi_{i}=\phi_{i}\circ f_{i+1}.$ Moreover because
of the truncation the exact proof used is irrelevant. For this reason
we conclude that if we have a map $\|\colim_{i}(Y_{i},\psi_{i})\|_{0}\rightarrow\|\colim_{i}(Z_{i},\epsilon_{i})\|_{0}$
given by some $g_{i}:Y_{i}\rightarrow Z_{i}$ satisfying the appropriate
relations then the map $\|\colim_{i}(X_{i},\phi_{i})\|_{0}\rightarrow\|\colim_{i}(Z_{i},\epsilon_{i})\|_{0}$
given by $f_{i}\circ g_{i}$ is the same as we would get from composing
the two maps $\|\colim_{i}(X_{i},\phi_{i})\|_{0}\rightarrow\|\colim_{i}(Y,\psi_{i})\|_{0}\rightarrow\|\colim_{i}(Z_{i},\epsilon_{i})\|_{0}$
. Finally if the $f_{i}$ are isomorphisms then the maps we get are
isomorphisms.
\end{remark}
\begin{theorem}
$H_{n}(-)$ is a functor.
\end{theorem}
\begin{proof}
First let's be clear about the functorial action. Since $\pi_{n+i}^{s}(-\land K_{i})$
is a functor for each $i$ it is easy to define an action on functions.
We just need to check that the following commutes:

\begin{align*}
\susp(X\land K_{i})\circ\pi_{n+i+1}^{s}(\suspsmtwo(X,K_{i})\circ X\land\kappa_{i})\circ\pi_{n+i+1}^{s}(f\land K_{i+1})\\
=\pi_{n+i}^{s}(f\land K_{i})\circ\susp(Y\land K_{i})\circ\pi_{n+i+1}^{s}(\suspsmtwo(Y,K_{i})\circ Y\land\kappa_{i}).
\end{align*}
We use that $\susp$ is natural and then cancel after which it suffices
to show

\begin{align*}
\pi_{n+i+1}^{s}(\suspsmtwo(X,K_{i})\circ X\land\kappa_{i})\circ\pi_{n+i+1}^{s}(f\land K_{i+1})\\
=\pi_{n+i+1}^{s}(\Sigma(f\land K_{i}))\circ\pi_{n+i+1}^{s}(\suspsmtwo(Y,K_{I})\circ Y\land\kappa_{i}).
\end{align*}
By functoriality it suffices to show

\[
\suspsmtwo(X,K_{i})\circ X\land\kappa_{i}\circ f\land K_{i+1}=(\Sigma(f\land K_{i}))\circ\suspsmtwo(Y,K_{i})\circ Y\land\kappa_{i}.
\]
By Lemma \ref{lem:suspsmash} 
\[
X\land\kappa_{i}\circ f\land K_{i+1}=f\land\Sigma K_{i}\circ Y\land\kappa_{i},
\]
and so it remains to show 
\[
\suspsmtwo(X,K_{i})\circ f\land\Sigma K_{i}=(\Sigma(f\land K_{i}))\circ\suspsmtwo(Y,K_{i}),
\]
which again follows by the same techniques mentioned in Lemma \ref{lem:suspsmash}.
The fact that this preserves composition is immediate from the previous
remark and the functoriality of $\pi_{n+i}^{s}(-\land K_{i})$.
\end{proof}
\begin{theorem}
$H_{n}(-)$ satisfies the suspension axiom.
\end{theorem}
\begin{proof}
We use $\susptwo$ on each component. By the above remark it suffices
to show the following commutes:

\begin{align*}
\susp(n+i,X\land K_{i})\circ\pi_{n+i+1}^{s}(\suspsmtwo(X,K_{i})\circ X\land\kappa_{i})\\
\circ\susptwo(n+1+i,X,K_{i+1})\\
=\susptwo(n+i,X,K_{i})\circ\susp(n+i+1,(\Sigma X)\land K_{i})\circ\\
\pi_{n+i+2}^{s}(\suspsmtwo(\Sigma X,K_{i})\circ\Sigma X\land\kappa_{i}).
\end{align*}
Expand the definition of $\susptwo$ and cancel the first term. We
get
\begin{align*}
\pi_{n+i+1}^{s}(\suspsmtwo(X,K_{i})\circ X\land\kappa_{i})\circ\susp(n+1+i,X\land K_{i+1})\circ\\
\pi_{n+i+2}(\suspsm(X,K_{i+1})\\
=\pi_{n+i+1}^{s}(\suspsm(X,K_{i}))\circ\susp((n+i+1,(\Sigma X)\land K_{i})\circ\\
\pi_{n+i+2}^{s}(\suspsmtwo(\Sigma X,K_{i})\circ\Sigma X\land\kappa_{i})
\end{align*}
We can then move the last term of the left hand side to the right
and use functoriality of $\pi_{n+i+2}$ to get $\pi_{n+i+2}^{s}(\suspsmtwo(\Sigma X,,K_{i})\circ\Sigma X\land\kappa_{i}\circ\suspsm(X,K_{i+1})^{-1})$
on the right hand side. Then we use Lemma \ref{lem:suspsmash} to
show
\begin{align*}
\suspsmtwo(\Sigma X,,K_{i})\circ\Sigma X\land\kappa_{i}\circ\suspsm(X,K_{i+1})^{-1}\\
=\Sigma\big(\suspsm(X,K_{i})^{-1}\circ\suspsmtwo(X,K_{i})\circ(X\land\kappa_{i})\big).
\end{align*}
So we can apply functoriality of $\Sigma$ and naturality of $\susp$
to get 
\begin{align*}
\pi_{n+i+1}^{s}(\suspsmtwo(X,K_{i})\circ X\land\kappa_{i})\circ\susp(n+1+i,X\land K_{i+1})
\end{align*}
is equal to 
\begin{align*}
\pi_{n+i+1}^{s}(\suspsm(X,K_{i}))\circ\pi_{n+i+1}^{s}(\suspsm(X,K_{i})^{-1}\\
\circ\suspsmtwo(X,K_{i})\circ(X\land\kappa_{i}))\circ\susp(n+i+1,X\land K_{i+1}).
\end{align*}
Finally we are done by simple canceling and using functoriality of
$\pi_{n+i+1}^{s}$.

Naturality follows from the above remark and the fact that $\susptwo$
is natural.

We conclude by noting exactness follows easily from exactness of $\pi_{n}^{s}(-\land K)$
and basic arguments about truncated colimits. So in conclusion we
have: 
\end{proof}
\begin{corollary}
\label{cor:colim-is-homo}$H_{n}(-)$ is a homology theory.
\end{corollary}
\begin{remark}
We get regular homology with coefficients in $G$ by using the spectrum
$K_{i}=K(G,i)$, the Eilenberg-Maclane spaces. These have been already
defined in Homotopy type theory \cite{licata2014eilenberg}.
\end{remark}
\section{Closing Remark}

Notice that the proof for the following (infinite) additivity axiom remains open.\begin{question}For a set $I$ and family $X:I \rightarrow \mathcal{U}_\cdot $ is the canonical map
\[\bigoplus_i H_n ( X_i )\rightarrow H_n(\bigvee_i X_i)\] an isomorphism? \end{question} To this end, it remains to show the smash product distributes over an infinite wedge product and that the sphere is compact. Of note is how we will deal with colimits of path spaces, a problem that has arisen in various contexts. See for example the issue of constructing spectra from prespectra \cite{cavallo2015synthetic} or the problem of computing the fundamental group of an infinite wedge.

 \begin{question}
 Given a family of pointed types $X:\mathbb{N}$
 $\rightarrow$
 $ \mathcal{U}_t$, does the following hold?
 \[ \Omega \colim_n  X(n)\]
 \[ \simeq  \colim_n \Omega  X(n) \]
 \end{question}


\bibliographystyle{plainnat}

\bibliography{Homology_in_HoTT}

\end{document}